\documentclass[12pt,twoside,leqno]{article}
\usepackage{amsmath}
\usepackage{amssymb}
\usepackage{amsxtra}
\usepackage{amscd}
\usepackage{amsthm}
\usepackage[mathscr]{eucal}

\theoremstyle{plain}
\newtheorem{thm}[subsection]{Theorem}
\newtheorem{prop}[subsection]{Proposition}
\newtheorem{cor}[subsection]{Corollary}
\newtheorem{lem}[subsection]{Lemma}

\theoremstyle{definition}

\newtheorem{rem}[subsection]{Remark}
\newtheorem{para}[subsection]{}

\newenvironment{pf}{\proof[\proofname]}{\endproof}

\begin{document}
\title{Logarithmic Dieudonn\'e theory}
\author
{Kazuya Kato}

\maketitle



\newcommand\Cal{\mathcal}
\newcommand\define{\newcommand}
\renewcommand\bold{\Bbb}

\define\bG{\bold G}
\define\bZ{\bold Z}
\define\bC{\bold C}
\define\bR{\bold R}
\define\bQ{\bold Q}
\define\bN{\bold N}
\define\bP{\bold P}
\define\cl{\mathrm{cl}}%
\define\fil{\mathrm{fil}}%
\define\fl{\mathrm{fl}}%
\define\gp{\mathrm{gp}}%
\define\fs{\mathrm{fs}}%
\define\an{\mathrm{an}}%
\define\et{\mathrm{\acute{e}t}}
\define\mult{\mathrm{mult}}%
\define\sat{\mathrm{sat}}%
\renewcommand\int{\mathrm{int}}%
\define\Ker{\mathrm{Ker}\,}%
\define\Coker{\mathrm{Coker}\,}%
\define\Hom{\operatorname{\mathrm{Hom}}}%
\define\Aut{\operatorname{\mathrm{Aut}}}%
\define\Mor{\operatorname{\mathrm{Mor}}}%
\define\rank{\mathrm{rank}\,}%
\define\gr{\mathrm{gr}}%
\define{\cHom}{\operatorname{\mathcal{H}\mathit{om}}}
\define{\HOM}{\cHom}
\define{\cExt}{\operatorname{\mathcal{E}\mathit{xt}}}
\define{\cMor}{\operatorname{\mathcal{M}\mathit{or}}}
\define\cO{\Cal O}
\renewcommand{\O}{\cO}
\define{\cA}{\Cal A}
\define{\cC}{\Cal C}
\define{\cD}{\Cal D}
\define\cS{\Cal S}
\define\cL{\Cal L}
\define\cM{\Cal M}
\define\cG{\Cal G}
\define\cH{\Cal H}
\define\cE{\Cal E}
\define\cF{\Cal F}

\newcommand{\N}{{\mathbb{N}}}
\newcommand{\Q}{{\mathbb{Q}}}
\newcommand{\Z}{{\mathbb{Z}}}
\newcommand{\R}{{\mathbb{R}}}
\newcommand{\C}{{\mathbb{C}}}
\newcommand{\F}{{\mathbb{F}}}
\newcommand{\ol}[1]{\overline{#1}}
\newcommand{\respect}{\rightsquigarrow}
\newcommand{\compatible}{\leftrightsquigarrow}
\newcommand{\upc}[1]{\overset {\lower 0.3ex \hbox{${\;}_{\circ}$}}{#1}}
\newcommand{\Gmlog}{\bG_{m, \log}}
\newcommand{\Gm}{\bG_m}
\newcommand{\Ga}{\bG_a}
\newcommand{\ep}{\varepsilon}
\newcommand{\pe}{\frak p}
\newcommand{\Spec}{\operatorname{Spec}}
\newcommand{\val}{{\mathrm{val}}}
\newcommand{\lcf}{{\mathrm{lcf}}}
\newcommand{\Et}{{\mathrm{Et}}}
\renewcommand{\b}{\langle\,,\,\rangle}
\newcommand{\n}{\operatorname{naive}}
\newcommand{\bs}{\operatorname{\backslash}}
\newcommand{\Lie}{\operatorname{Lie}}
\newcommand{\fin}{\operatorname{fin}}
\newcommand{\nar}{{\mathrm{nar}}}
\newcommand{\crys}{{\mathrm{crys}}}
\newcommand{\nilp}{{\mathrm{nilp}}}
\newcommand{\lff}{{\mathrm{lff}}}
\newcommand{\Gal}{{\mathrm{Gal}}}
\renewcommand{\emph}{\it}
\newcommand{\Map}{\operatorname{Map}}
\newcommand{\la}{\operatorname{\lambda}}
\newcommand{\La}{\operatorname{\Lambda}}

\begin{abstract} 
  We develop the theory of logarithmic $p$-divisible groups and the theory of logarithmic finite locally  free commutative group schemes. 
  \end{abstract}

\renewcommand{\thefootnote}{\fnsymbol{footnote}}
\footnote[0]{MSC2020: Primary 14A21; 
Secondary 14L05, 14L15}

In this paper, we develop the logarithmic versions of the theory of finite locally  free commutative group schemes and the theory of $p$-divisible groups.

This paper was written in 1992 but  some proofs were omitted. It was circulated as an incomplete preprint for a long time. Here we supply the omitted proofs and correct mistakes except that the proof of Thm. \ref{6.3} will be given in a forthcoming paper \cite{K3}. Since that preprint (quoted in this paper as the old version) was already referred to in several published works, we tried to preserve the original form and in particular to preserve the numberings of definitions and propositions. Some new statements (for example, \ref{7.3}) which are not contained in the old version but used to supply the proofs,  are added.

The author wishes to express his hearty thanks to T. Kajiwara, C. Nakayama and H. Zhao  who helped him  in the completion of this paper. 

The author is partially supported by NSF Award 2001182.

\section{Finite logarithmic group objects ; definitions}

 Throughout this paper, let $X$ be an fs log scheme (\cite{K2} \S1) whose underlying scheme is locally Noetherian.

Let $X^{\log}_{\fl}$ (resp. $X^{\log}_{\et}$, resp. $X^{\cl}_{\fl}$) be the category of fs log schemes over $X$ regarded as a site with the logarithmic flat (resp. logarithmic \'etale, resp. classical flat) topology in the sense of \cite{K2} \S2. (These log flat topology and log \'etale topology are also called Kummer log flat topology and Kummer log \'etale topology, respectively.)
Thus these sites are the same as categories, but the topologies are different.

\begin{para}\label{1.1} We denote by $(\fin/X)_c$ the category of finite locally free commutative group schemes over the underlying scheme of $X$. 
(The notation "c'' means "\underline{c}lassical objects''.) By endowing the objects of $(\fin/X)_c$ with the inverse images of the log structure of $X$, and by considering the sheaf on $X^{\log}_{\fl}$ (\cite{K2} 3.1) represented by them, we regard $(\fin/X)_c$ as a full subcategory of the category of sheaves of abelian groups on $X^{\log}_{\fl}$. 
\end{para} 

\begin{para}\label{1.2} We define four categories
$$\ref{1.2}.1 \quad \quad (\fin/X)_d\subset (\fin/X)_r\subset (\fin/X)_e\subset (\fin/X)_f$$
which are full subcategories of the category of sheaves of abelian groups on $X^{\log}_{\fl}$. (The letters $d,r,e,f$ 
are taken from the properties of the categories or from the properties of the objects of the categories; "is stable under the Cartier \underline{d}uality'', 
"is \underline{r}epresentable'', 
"is classical log \underline{e}tale locally'', "is classical log \underline{f}lat locally'', respectively.) These categories are enlargements of $(\fin/X)_c$. 

The definitions of  $(\fin/X)_r$, $(\fin/X)_e$, $(\fin/X)_f$ are as follows. (The definition of $(\fin/X)_d$ will be given slightly later.)
\end{para}

\begin{para}\label{1.3} {\bf Definition.} Let $G$ be a sheaf of abelian groups on $X^{\log}_{\fl}$.

(1) $G$ belongs to $(\fin/X)_f$ if there is a covering $Y\to X$ in $X^{\log}_{\fl}$ (\cite{K2} 2.3) such that the pullback of $G$ to $Y^{\log}_{\fl}$ belongs to $(\fin/Y)_c$. 

(2) $G$ belongs to $(\fin/X)_r$ if it belongs to $(\fin/X)_f$ and it is represented by an fs log scheme over $X$.

(3) $G$ belongs to $(\fin/X)_e$ if there is a covering $Y\to X$ in $X^{\log}_{\et}$  such that the pullback of $G$ to $Y^{\log}_{\fl}$ belongs to $(\fin/Y)_r$. 

\end{para}

\begin{prop}\label{1.4}  Let $G$ be a sheaf of abelian groups on $X^{\log}_{\fl}$. Then $G$ belongs to $(\fin/X)_r$ if and only if $G$ is represented by an fs log scheme over $X$ which is log flat, finite and of Kummer type.

\end{prop}

This follows from \cite{K2} 2.7(1) and from the log flat decent theory for properties of morphisms (\cite{K2} \S7). 

\begin{para}\label{1.5} To define $(\fin/X)_d$, we need some preparation on Cartier duality. 

For an object $G$ of $(\fin/X)_f$, the sheaf $\cH om(G, {\bf G}_m)$ on $X^{\log}_{\fl}$ also belongs to $(\fin/X)_f$ and
$$G\overset{\cong}\to \cH om(\cH om(G, {\bf G}_m), {\bf G}_m).$$
Indeed, working locally for the log flat topology, these are reduced to the classical Cartier duality. As we will see later in \ref{1.8}.2, $(\fin/X)_r$ is not stable under Cartier duality. 

\end{para}

\begin{para}\label{1.6} We define $(\fin/X)_d$ to be the full subcategory of $(\fin/X)_r$ consisting of objects $G$ such that $\cH om(G, {\bf G}_m)$ also belongs to $(\fin/X)_r$.

\end{para}

\begin{rem}\label{1.7} We will see in \ref{2.5} that $(\fin/X)_e$ is stable under $\cH om(G, {\bf G}_m)$. 

\end{rem}

\begin{para}\label{1.8} Now we give examples concerning the differences between the categories $(\fin/X)_c$, $(\fin/X)_d$, $(\fin/X)_r$, $\dots$. The proofs of the statements in this \ref{1.8} (which were not given in the old version) will be given in  \ref{pf1.8} after we prove in \S2 the basic results on these categories. 

\ref{1.8}.1. Let $n\geq 1$, and let $\delta_n: H^0(X, M_X^{\gp}) \to H^1(X^{\log}_{\fl}, \Z/n(1))$ be the connecting map of the Kummer sequence (\cite{K2} 4.2)
$$0\to \Z/n(1) \to {\bf G}_{m,\log}\overset{n}\to {\bf G}_{m,\log}\to 0.$$
Here ${\bf G}_{m,\log}$ is the sheaf $Y\mapsto \Gamma(Y, M_Y^{\gp})$ on $X^{\log}_{\fl}$ (\cite{K2}  3.2). 

Let $a\in \Gamma(X,  M_X^{\gp})$ and let 
$$0\to \Z/n(1) \to G \to \Z/n\to 0$$
be the extension corresponding to $\delta_n(a)\in H^1(X^{\log}_{\fl}, 
\Z/n(1))=\text{Ext}^1_{\Z/n}(\Z/n, \Z/n(1))$. Then $G$ belongs to $(\fin/X)_d$. If there exists a point $x\in X$ such that the image of $a$ in $(M_X^{\gp}/\cO_X^\times)_{\bar x}$ is not an $n$-th power, $G$ does not belong to $(\fin/X)_c$. 

\ref{1.8}.2. Let $p$ be a prime number. Assume that $p-1$ is invertible on $X$ and that a primitive $p-1$ th root of $1$ exists in $\Gamma(X, \cO_X)$. Then, the image of $\delta_{p-1}(a)$ under
$$H^1(X^{\log}_{\fl}, (\Z/(p-1))(1)) \cong H^1(X^{\log}_{\fl}, (\Z/p)^\times) \cong H^1(X^{\log}_{\fl}, {\cal A}ut(\Z/p))$$
defines a ''twist'' of $\Z/p$ on $X^{\log}_{\fl}$, i.e., a sheaf $G$ of abelian groups on $X^{\log}_{\fl}$ which is log flat locally isomorphic to $\Z/p$. Then, $G$ belongs to $(\fin/X)_r$ and its Cartier dual $\cH om(G, {\bf G}_m)$ 
belongs to $(\fin/X)_e$. Assume that there exists a point $x\in X$ such that the image of $a$ in 
$(M_X^{\gp}/\cO^\times_X)_{\bar x}$ is not a $(p-1)$-th power. Then $G$ does not belong to $(\fin/X)_d$. If furthermore $p$ is not invertible at this $x$,  $\cH om(G, {\bf G}_m)$ does not belong to $(\fin/X)_r$. 

\ref{1.8}.3. Let $p$ be a prime number and assume $p\cO_X=0$. Let $\cL$ be a line bundle on $X^{\log}_{\fl}$ and let $\cL^\times \subset \cL$ be the ${\bf G}_m$-torsor consisting of bases of $\cL$.  
Let  $\alpha_p=\text{Ker}({\bf G}_a\to {\bf G}_a\;;\; x\mapsto x^p)$. Since $\cA ut(\alpha_p)= {\bf G}_m$, $\cL^\times$ defines a twist $G$ of $\alpha_p$ in $(\fin/X)_f$.
The following (i), (ii) and (iii) are equivalent.

(i) $G$ belongs to $(\fin/X)_c$.

(ii) $G$ belongs to $(\fin/X)_r$.

(iii) $\cL$ comes from a line bundle on $X_{\text{Zar}}$.

The following (iv) and (v) are equivalent.

(iv) $G$ belongs to $(\fin/X)_e$.

(v) For every $x\in X$, the class of $\cL$ in 
 $\Q/\Z \otimes (M_X^{\gp}/\cO_X^\times)_{\bar x}$ (\cite{K2} 6.2) at the strict henselization at $x$ has order coprime to $p$. 
 
\end{para}

\section{Finite logarithmic group objects ; basic results}

If $G$ is an object of $(\fin/X)_f$ and $X$ is quasi-compact, $G$ is killed by some non-zero integer. In fact, since there is a covering $Y\to X$ for the log flat topology such that $Y$ is quasi-compact and the pullback of $G$ to $Y$ is a classical finite locally free commutative group scheme, this is reduced to the classical result on this pullback.

The proofs of the following results \ref{2.1}, \ref{2.3}, \ref{2.4} will be given in \ref{pf2.1}, \ref{pf2.3}, \ref{pf2.4}, respectively, after local studies.

\begin{prop}\label{2.1} Let $G$ be an object of $(\fin/X)_f$ and assume that $G$ is killed by an integer which is invertible on $X$. Then $G$ belongs to $(\fin/X)_d$. 

\end{prop}

\begin{cor}\label{2.2} If the underlying scheme of $X$ is a $\Q$-scheme, the categories $(\fin/X)_d$, $(\fin/X)_r$, $(\fin/X)_e$ and $(\fin/X)_f$ coincide.

\end{cor}

\begin{prop}\label{2.3}

Let $0\to G'\to G \to G'' \to 0$ be an exact sequence of sheaves of abelian 
groups on $X^{\log}_{\fl}$. If $G'$ and $G''$ belong to $(\fin/X)_d$ (resp. $(\fin/X)_r$, resp. $(\fin/X)_e$, resp. $(\fin/X)_f$), so does $G$. 
(The similar statement for $(\fin/X)_c$ is not true, and this matter is studied in \S3.)
\end{prop}

\begin{prop}\label{2.4} Let $G$ be a sheaf of abelian groups on $X^{\log}_{\fl}$. Then the following three conditions (i) - (iii) are equivalent.

(i) $G$ belongs to $(\fin/X)_e$.

(ii) Log \'etale locally on $X$, there exists a sequence $0=H_0\subset H_1 \subset \dots \subset H_r=G$ of subgroup sheaves of $G$ for some $r\geq 0$ such that $H_i/H_{i-1}$ belongs to $(\fin/X)_c$ for $1\leq i\leq r$.

(iii) Classical \'etale locally on $X$, there exists a subgroup sheaf $H$ of $G$ such that log \'etale locally, $H$ and $G/H$ belong to $(\fin/X)_c$.

\end{prop}

\begin{cor}\label{2.5} $(\fin/X)_e$ is stable under the Cartier duality.

\end{cor}

\begin{para}\label{2.6} We consider a local situation. Assume that as a scheme, $X=\Spec(A)$ for a Noetherian strict local ring, and let $G$ be an object of $(\fin/X)_f$. Then
there exists a unique exact sequence of objects of $(\fin/X)_f$
$$0 \to G^{\circ} \to G \to G^{\et} \to 0$$
having the following property. There exists a covering $Y\to X$ in $X^{\log}_{\fl}$ such that $Y$ is finite over $X$ and connected (hence $Y$ is also the Spec of a Noetherian strict local ring), and such that the pullback of $G^{\circ}$ (resp. $G^{\et}$) on $Y^{\log}_{\fl}$ is represented by a connected (resp. \'etale) finite locally free commutative group scheme over the underlying scheme of $Y$, which is endowed with the inverse image of $M_Y$.

Indeed, take a chart $P\to M_X$. Then there exists an fs monoid $Q$ and a homomorphism $P\to Q$ of Kummer type  such that, if we denote $X\otimes_{\Z[P]} \Z[Q]$ by $Y$ and endow it with the log structure associated to $Q\to \cO_Y$, then the pullback of $G$ on $Y^{\log}_{\fl}$ belongs to $(\fin/Y)_c$ (\cite{K2} 2.7).  The usual connected-\'etale exact sequence for $G\times_X Y$ descends to an exact sequence of sheaves on $X^{\log}_{\fl}$ as is seen by working on $Y\times_X Y$.

\end{para}

\begin{prop}\label{2.7} Assume $X=\Spec(A)$ as a scheme for a Noetherian strict local ring $A$. Let $G$ be an object of $(\fin/X)_f$ and let $G^{\circ}$ and $G^{\et}$ be as in \ref{2.6}. 

(1) $G^{\et}$ belongs to $(\fin/X)_r$, and is log \'etale over $X$.

(2) $G$ belongs to $(\fin/X)_r$ if and only if $G^{\circ}$ belongs to $(\fin/X)_c$. If $G$ belongs to 
$(\fin/X)_r$, $G^{\circ}$ is the connected component of the scheme $G$ containing the origin of $G$ endowed with the inverse image of the log structure of $G$.

(3) Assume that the characteristic $p$ of the residue field of $A$ is non-zero and $G$ is killed by some power of $p$. Then, $G$ belongs to $(\fin/X)_d$ if and only if $G^{\circ}$ and $G^{\et}$ belong to $(\fin/X)_c$.

\end{prop}

\begin{pf} Here we prove (1) and the only if part of (2). The proof of the if part of (2), and that of (3) will be given  in \ref{pf2.7.2if} and \ref{pf2.7.3}, respectively.

(1) follows from \cite{K2} 10.2.

 We prove the only if part of (2).

Assume $G$ belongs to $(\fin/X)_r$. Then $G^{\circ}$ is representable since $G^{\circ}$ is the kernel of the homomorphism $G\to G^{\et}$ of representable objects (\ref{2.7} (1)) 
and since the category of fs log schemes has fiber products. 
We next prove $G^{\circ}$ is connected. 
Take a covering $Y\to X$ in $X^{\log}_{\fl}$  as in 2.6.
 Then $G^{\circ}\times_X Y$ is connected and $G^{\circ}\times_X Y \to G^{\circ}$ is surjective, and hence $G^{\circ}$ is connected. Now, by Lemma \ref{2.8} below, $G^{\circ}$ belongs to $(\fin/X)_r$. 

Let $U$ be the connected component of $G$ containing  the origin in $G$,  endowed with the inverse image of the log structure of $G$. Since $G^{\circ}$ is connected, the canonical 
homomorphism $G^{\circ}\to G$ factors through $U$. We show $G^{\circ}\overset{\cong}\to U$. For $Y$ as above, we have $G^{\circ}\times_X Y \overset{\cong}\to U \times_X Y$. 
This implies $G^{\circ}\overset{\cong}\to U$ since $G^{\circ}$ 
is the quotient of $G^{\circ} \times_X Y$ by the 
equivalence relation given by the two projections from $G^{\circ} \times_X Y \times_X Y$ to $G^{\circ} \times_X Y$ and since we have the similar understanding of $U$. 
\end{pf}

\begin{lem}\label{2.8} Let $X$ be as in the hypothesis of \ref{2.7}, and let $G$ be an object of $(\fin/X)_r$. Assume $G$ is connected. Then, $G$ belongs to $(\fin/X)_c$.

\end{lem}
\begin{pf} It is sufficient to show that the log structure of $G$ is the inverse image of that of $X$. Let $x$ be the closed point of $X$. The origin $X \to G$ of $G$ defines a homomorphism $(M_G/\cO_G^\times)_y \to (M_X/\cO_X^\times)_x$, where $y$ is the image of $x$ in $G$.  Since the composite $(M_X/\cO_X^\times)_x \to (M_G/\cO_G^\times)_y \to 
(M_X/\cO_X^\times)_x$ is the identity map and the first arrow is a homomorphism of fs monoids of Kummer type, we have $(M_X/\cO_X^\times)_x \overset{\cong}\to (M_G/\cO_G^\times)_y$. Since $y$ is the only closed point in $G$ (for $G$ is connected), this implies that $M_G$ is the inverse image of $M_X$.
\end{pf}

\begin{lem}\label{2.9} Let $0\to G'\to G \to G'' \to 0$ be an exact sequence of sheaves of abelian 
groups on $X^{\log}_{\fl}$, and assume that $G'$ belongs to $(\fin/X)_c$ and $G'' $ belongs to $(\fin/X)_r$. Then $G$ belongs to $(\fin/X)_r$. 

\end{lem}

\begin{pf} Since $G$ is a $G'$-torsor over $G''$, this follows from \cite{K2} 9.1. 
\end{pf}

\begin{para}\label{pf2.7.2if} The if part of \ref{2.7} (2) follows from \ref{2.7} (1) and Lemma \ref{2.9}.
\end{para}

\begin{para}\label{pf2.3} We prove \ref{2.3}. It is sufficient to consider $(\fin/X)_r$. But we first notice that the statement for $(\fin/X)_f$ follows from \ref{2.9} by a log flat local argument.

To prove the statement for $(\fin/X)_r$, by the classical \'etale descent and by a limit argument, we may assume that $X=\Spec(A)$ as a scheme for a Noetherian strict local ring $A$. Assume $G'$ and $G''$ belong to $(\fin/X)_r$. Then $G$ belongs to $(\fin/X)_f$ and we have an exact sequence $0\to (G')^{\circ}\to G^{\circ} \to (G'')^{\circ}\to 0$. By \ref{2.7} (2), $(G')^{\circ}$ and $(G'')^{\circ}$ belong to $(\fin/X)_c$ and our task is to prove that $G^{\circ}$ belongs to $(\fin/X)_c$. Let $\epsilon: X^{\log}_{\fl}\to X^{\cl}_{\fl}$ be the canonical morphism of sites. Then we have an exact sequence on $X^{\cl}_{\fl}$
$$0\to (G')^{\circ} \to G^{\circ} \to (G'')^{\circ}\overset{\delta}\to \cH$$
where $$\cH=R^1\epsilon_*((G')^{\circ})\cong \varinjlim_n  \cH om(\Z/n(1), (G')^{\circ})\otimes {\bf G}_{m,\log}/{\bf G}_m$$  (\cite{K2} 4.1). 

We prove $\delta=0$. 

 Let $Y:=(G'')^{\circ}$. Then $Y$ is the Spec of a Noetherian strict local ring. Let $x$ be the closed point of $X$ regarded as the Spec of the residue field of $A$ endowed with the inverse image of the log structure of $X$. 
  We have a commutative diagram
$$\begin{matrix} \Hom((G'')^{\circ}, \cH) &\overset{\subset}\to  & \Gamma(Y, \cH)\\
\downarrow && \downarrow \;\cong\\
\Hom((G'')^{\circ}|_x, \cH|_x) &\overset{\subset}\to  & \Gamma(Y \times_X x, \cH),
\end{matrix}$$
where $|_x$ denotes the pullback to $x^{\cl}_{\fl}$. The right vertical arrow is an isomorphism because, if we denote the multiplicative part of $(G')^{\circ}$
by $(G')^{\mult}$ and take  $n\geq 1$ which kills $(G')^{\circ}$, $\cG:=\cHom(\Z/n(1), (G')^{\mult})$ is an \'etale finite group scheme over $X$, and we have
$\Gamma(Y,\cH)= \Gamma(Y, \cG\otimes {\bf G}_{m,\log}/{\bf G}_m) \cong \Gamma(Y \times_X x, \cG\otimes {\bf G}_{m,\log}/{\bf G}_m)= \Gamma(Y\times_X x, \cH)$.

By this diagram, the map $\Hom((G'')^{\circ}, \cH)\to \Hom((G'')^{\circ}|_x, \cH|_x)$ is injective. But $\cH|_x$ is the constant sheaf associated to the group
$\cH om(\Z/n(1), (G')^{\circ}|_x)\otimes (M_X^{\gp}/\cO_X^\times)_x$  Since $(G'')^{\circ}|_x$ is connected, we have 
$\Hom((G'')^{\circ}|_x, \cH|_x)=0$.

 Since $\delta=0$, we have 
an exact sequence $0\to (G')^{\circ} \to G^{\circ}\to (G')^{\circ} \to 0$ on $X^{\cl}_{\fl}$. This shows that $G^{\circ}$ belongs to $(\fin/X)_c$. 

\end{para}

\begin{para}\label{pf2.1} We prove \ref{2.1}. We may assume that as a scheme, $X=\Spec(A)$ for a Noetherian strict local ring $A$. By the assumption, we have $G^{\circ}=1$. Hence by \ref{2.7} (2), $G$ belongs to $(\fin/X)_r$. By considering the Cartier dual, we see that $G$ belongs to $(\fin/X)_d$.

\end{para}

\begin{para}\label{pf2.4} We prove \ref{2.4}. We prove that (i) implies (iii).  We may assume that $X=\Spec(A)$ for a Noetherian strict local ring $A$. Then take $H=G^{\circ}$. That (iii) implies (ii) is clear. Finally (ii) implies (i) by \ref{2.3}.

\end{para}

\begin{para}\label{pf2.7.3} We prove \ref{2.7} (3). Let $G$ be an object of $(\fin/X)_d$ which is killed by some power of $p$. We know already that $G^{\circ}$ belongs to $(\fin/X)_c$ (\ref{2.7} (2)), and so we prove that $G^{\et}$ belongs to $(\fin/X)_c$. Let $H$ be the Cartier dual of $G$. Then, $G^{et}$ is the Cartier dual of the multiplicative part $H^{\text{mult}}$ of $H^{\circ}$. Since $H^{\circ}$ belongs to $(\fin/X)_c$ (\ref{2.7} (2)), $H^{\text{mult}}$ also belongs to $(\fin/X)_c$ and hence its Cartier dual $G^{\et}$ belongs to $(\fin/X)_c$.

Conversely, if $G^{\circ}$ and $G^{\et}$ belong to $(\fin/X)_c$, then $G$ belongs to $(\fin/X)_d$ by \ref{2.3}. 
\end{para}

\begin{prop}\label{2OG}   Let $G$ be an object of $(\fin/X)_f$. Define the sheaf $\cO_G$ on $X^{\log}_{\fl}$ by
$$\cO_G= \cM or(G, \cO_X),$$
where $\cO_X$ denotes the sheaf $T\mapsto \Gamma(T, \cO_T)$ on $X^{\log}_{\fl}$ and $\cM or$ denotes the sheaf of morphisms as sheaves of sets (the group structures are forgotten here). Then:

(1) As an $\cO_X$-module on $X^{\log}_{\fl}$, $\cO_G$ is locally free of finite rank. 

(2)  $G$ belongs to $(\fin/X)_c$ if and only if the $\cO_X$-module $\cO_G$ is classical (that is, locally free of finite rank for the classical Zariski topology of $X$). 
\end{prop}

\begin{pf} 

(1) By log flat localization, we are reduced to the classical case.

(2) The only if part is clear. We prove the if part. By the reduction to the classical case, we see that $G$ is the sheaf of homomorphisms $\cO_G \to \cO_X$ of $\cO_X$-algebras. Hence if $\cO_G$ comes from a sheaf $R$ of $\cO_X$-algebras on the Zariski site $X_{\text{zar}}$ which is locally free of finite type as an $\cO_X$-module, then $G$ is represented by $\Spec(R)$ with the inverse image of the log structure of $X$ and hence belongs to $(\fin/X)_c$.
\end{pf}

\begin{prop}\label{2pt}  Let $G$ be an object of $(\fin/X)_f$. Assume that for every closed point $x$ of $X$, the pullback of $X$ on $x$ ($=$ the Spec of the residue field of $x$ endowed with the inverse image of $M_X$)
belongs to $(\fin/x)_c$ (resp. $(\fin/x)_d$, resp. $(\fin/x)_r$, resp. $(\fin/x)_e$). Then $G$ belongs to $(\fin/X)_c$ (resp. $(\fin/X)_d$, resp. $(\fin/X)_r$, resp. $(\fin/X)_e$).

\end{prop}

\begin{pf}  
We prove the statement for $(\fin/X)_c$. By \ref{2OG}, it is enough to prove that $\cO_G$ is classical. We may assume that the underlying scheme of $X$ is the Spec of a Noetherian strict local ring. Then we are reduced to 
 \cite{K2} Thm. 6.2.

We prove the  statements for $(\fin/X)_r$ and for $(\fin/X)_e$. We may assume that the underlying scheme of $X$ is the Spec of a Noetherian strict local ring. By \ref{2.7} (2), the statement for $(\fin/X)_r$ in this case is reduced to the statement for $(\fin/X)_c$ applied to $G^{\circ}$. The statement for $(\fin/X)_e$ in this case is reduced to the statement for $(\fin/X)_r$.

The statement for $(\cdot)_d$ follows from that for $(\fin/X)_r$.
\end{pf}

\begin{para}\label{pf1.8} We prove the statements in \ref{1.8}.

\ref{1.8}.1.  Since $G$ is a $\Z/n(1)$-torsor over $\Z/n$, it belongs to $(\fin/X)_r$ by  \cite{K2} 9.1. By the isomorphism $\wedge^2_{\Z/n}\; G \cong \Z/n(1)$, the Cartier dual of $G$ is isomorphic to $G$. Hence $G$ belongs to $(\fin/X)_d$. 

Assume $G$ belongs to $(\fin/X)_c$. Then the inverse image 
$T\subset G$  of $1\in \Z/n$ satisfies $\Z/n(1) \times T \overset{\cong}\to  T \times_X T \;;\; (a,b) \mapsto (b, ab)$ and $T\to X$ is surjective, and hence $T$ is a $\Z/n(1)$-torsor in the classical sense. Hence if $\epsilon: X^{\log}_{\fl}\to X^{\cl}_{\fl}$ denotes the morphism of sites, the image of $\delta_n(a)$ in $R^1\epsilon_*(\Z/n(1))\cong \Z/n \otimes M^{\gp}_X/\cO_X^\times$ is zero. This shows that
 for each $x\in X$, 
 the image of $a$ in $(M_X^{\gp}/\cO_X^\times)_{\bar x}$ is  an $n$-th power.

\ref{1.8}.2. Since $p-1$ is invertible on $X$, the Kummer sequence $0\to \Z/(p-1) \to {\bf G}_{m,\log}\overset{p-1}\to {\bf G}_{m,\log}\to 0$ on $X^{\log}_{\et}$ is exact and hence $\delta_{p-1}(a)$ comes from $H^1(X^{\log}_{\et}, \Z/(p-1))$. Hence $\delta_{p-1}(a)$ vanishes log \'etale locally on $X$, and hence $G$ is log \'etale locally classical. Hence $G$ and its Cartier dual belong to $(\fin/X)_e$. Let $T$ be the $(\Z/p)^\times$-torsor representing $\delta_{p-1}(a)$. Then $G=T\times^{(\Z/p)^\times} \Z/p= T \times^{(\Z/p)^\times} (\{0\} \coprod (\Z/p)^\times)  \cong X \coprod T$. Since $T$ is representable by \cite{K2} 9.1, $G$ belongs to $(\fin/X)_r$. By \ref{2.5}, the Cartier dual $G^*$ of $G$ belongs to $(\fin/X)_e$.

Assume that $G^*$ belongs to $(\fin/X)_r$. Then since $G^*=(G^*)^{\circ}$, $G$ belongs to $(\fin/X)_c$ by \ref{2.7} (2). Hence the Cartier dual $G$ of $G^*$ belongs to $(\fin/X)_c$. This proves that the log structure of the above $T$ is the inverse image of that of $X$. Hence $T$ is a classical $\Z/(p-1)$-torsor and hence $\delta_{p-1}(a)$ in $R^1\epsilon_*(\Z/(p-1))\cong \Z/(p-1) \otimes M_X^{\gp}/\cO_X^\times$ is zero. Hence 
 for each $x\in X$, 
 the image of $a$ in $(M_X^{\gp}\cO_X^\times)_{\bar x}$ is  an $p-1$-th power.

\ref{1.8}.3. We have $G=(\cL^\times)\times ^{{\bG}_m} \alpha_p$. Hence $$\cO_G\cong \cM or((\cL^\times)\times^{{\bG}_m} \alpha_p,\; \cO_X) =  \oplus_{i=0}^{p-1} \; \cL^{\otimes -i},$$ where  $(\ell_i)_{0\leq i\leq p-1}$ on the right hand side corresponds to the morphism $(e, a) \mapsto \sum_i \ell_i(e^{\otimes i})a^i$ ($e\in \cL^\times, a\in \alpha_p)$. Hence $\cO_G$ is a classical vector bundle if and only if $\cL$ is a classical line bundle. Hence by \ref{2OG}, the conditions (i) and (iii) are equivalent. By \ref{2.7} (2), since $G=G^{\circ}$, the conditions (i) and (ii) are equivalent.

By this, $G$ belongs to $(\fin/X)_e$ if and only if $\cL$ comes from a line bundle on $X^{\log}_{\et}$. Hence by \cite{K2} 6.2, the conditions (iv) and (v) are equivalent.

\end{para}

\section{Finite logarithmic group objects ; a local equivalence}
In this \S3, we consider the category $(\fin/X)_d$.

\begin{thm}\label{3.1}
Assume $X=\Spec(A)$ as a scheme for a Noetherian strict local ring $A$. Assume that the characteristic $p$ of the residue field of $A$ is non-zero. Fix a chart $P\to M_X$ such that $P\overset{\cong}\to (M_X/\cO_X^\times)_x$ where $x$ is the closed point of $X$. Then there exists an equivalence of categories 
$$\cC \simeq \cC',$$
where $\cC$ is the full subcategory of $(\fin/X)_d$ consisting of objects which are killed by some powers of $p$, and $\cC'$ is the category of pairs $(G, N)$ with $G$ an object of $(\fin/X)_c$ killed by some power of $p$ and $N$ is a homomorphism
$$G^{\et}(1) \to G^{\circ} \otimes P^{\gp}.$$

\end{thm}

Cf. \ref{2.6} for $G^{\circ}$ and $G^{\et}$. Here $G^{\et}(1)= G^{\et}\otimes_{\Z/n} \Z/n(1)$ for a non-zero integer $n$ which kills $G^{\et}$. In other words, $G^{\et}(1)= \cH om(\cH om(G^{\et}, \Q/\Z), {\bf G}_m)$. Next, $G^{\circ} \otimes P^{\gp}$ is not a complicated thing. If we fix a $\Z$-basis $(e_i)_{1\leq i\leq r}$ of $P^{\gp}$, $G^{\circ} \otimes P^{\gp}$ is identified with the product $\prod^r G^{\circ}$ of $r$ copies of $G^{\circ}$ via $\prod^r G^{\circ} \overset{\cong}\to G^{\circ} \otimes P^{\gp}\;;\; (x_i)_i \mapsto \sum_i x_i \otimes e_i$. 

The functor between $\cC$ and $\cC'$ defined below, which gives the equivalence,  depends on the choice of the chart $P\to M_X$. 

The proof of Thm. \ref{3.1} was explained in the original version of this paper, but there was a mistake in \ref{3.7}. A corrected proof was given in 
\cite{WZ} \S3.1. Here for the convenience of the reader, we give the proof preserving the original form except that the mistake is  corrected as in \cite{WZ}. (A similar mistake existed in Section 2 of the old version, and is corrected in \ref{pf2.3} of this version.)

\begin{para}\label{3.2} We will deduce Thm. \ref{3.1} from Thm. \ref{3.3} below. We introduce some notation used in \ref{3.3}.

Let $X$ be an fs log scheme and let $G, H$ be objects of $(\fin/X)_c$. We denote by
$${\frak E}xt_{\log}(G, H)    \quad (\text{resp}.  \;\; {\frak E}xt_{\cl}(G, H))  $$
the category of sheaves of abelian groups $E$ on $X^{\log}_{\fl}$ (resp. $X^{\cl}_{\fl}$) endowed with an exact sequence $0\to H \to E \to G\to 0$. The set of isomorphism classes of ${\frak E}xt_{\log}(G, H)$ (resp. ${\frak E}xt_{\cl}(G, H)$) is $\text{Ext}^1(G, H)$ for the category of sheaves of abelian groups on $X^{\log}_{\fl}$ (resp. $X^{\cl}_{\fl}$). 

\end{para}

\begin{thm}\label{3.3} Let $X$, $A$ and $P$ be as in the hypothesis of \ref{3.1}. For objects $G$, $H$ of $(\fin/X)_c$, we have an equivalence of categories
$${\frak E}xt_{\log}(G, H) \simeq {\frak E}xt_{\cl}(G, H) \times (\Hom(G^{\et}(1), H)\otimes P^{\gp}).$$
\end{thm}
Here $\Hom(G^{\et}(1), H) \otimes P^{\gp}$ is regarded as a discrete category (an object is an element of this group and all morphisms are only the identity morphisms). The functor which gives the equivalence in \ref{3.3} defined below depends on the choice of the chart $P\to M_X$.

\begin{para}\label{3.4} We define the functor 

\medskip

\ref{3.4}.1.\;\; ${\frak E}xt_{\cl}(G, H) \times (\Hom(G^{\et}(1), H) \otimes P^{\gp}) \to {\frak E}xt_{\log}(G,H)$

\medskip
\noindent
by taking the Baer sum in ${\frak E}xt_{\log}(G, H)$ of the functors

\medskip

\ref{3.4}.2\;\; $\alpha: {\frak E}xt_{\cl}(G, H) \to {\frak E}xt_{\log}(G,H)$

\medskip

\ref{3.4}.3\;\;  $\beta: \Hom(G, H)\otimes P^{\gp} \to {\frak E}xt_{\log}(G,H)$

\medskip
\noindent
which are defined in \ref{3.5} and \ref{3.6} below, respectively.

\end{para}

\begin{para}\label{3.5} First we define the functor $\alpha$ by regarding an exact sequence $0\to H\to E \to G \to 0$ on $X^{\cl}_{\fl}$ as an exact sequence on $X^{\log}_{\fl}$.
Note that $G$ and $H$ are the functors represented by schemes over $A$ endowed with the inverse images of $M_X$ and so $G$ and $H$ are also sheaves on $X^{\log}_{\fl}$  by \cite{K2} 3.1. This functor 
$\alpha$ is clearly fully faithful.

\end{para}

\begin{para}\label{3.6} We define the functor $\beta$. For $h\in\Hom(G^{\et}(1), H)$ and $a\in P^{\gp}$, we define an object $\beta(h, a)$ of ${\frak E}xt_{\log}(G,H)$ as follows. 

Working locally on $X$, take an integer $n\geq 1$ which kills $G$. We define a sheaf of abelian groups $L$ on $X^{\log}_{\fl}$ by the commutative diagram of exact sequences 
$$\begin{matrix} 0& \to & \Z/n(1) & \to& {\bf G}_{m,\log} & \overset{n}\to & {\bf G}_{m,\log}& \to & 0\\
&& \parallel&& \uparrow && {\;}\uparrow{a} &&\\
0& \to & \Z/n(1) & \to& L& \overset{n}\to & \Z& \to & 0.
\end{matrix}$$
Let $E_a:=L/nL$. Then, we have an exact sequence $0\to \Z/n(1) \to E_a \to \Z/n \to 0$. By taking $\otimes_{\Z/n}\; G$, we obtain an exact sequence 
$0\to \Z/n(1) \otimes_{\Z/n} G \to E_a \otimes_{\Z/n} G \to G \to 0$. We define $\beta(h, a)$ as the pushout of $H \leftarrow \Z/n(1) \otimes_{\Z/n} G \to E_a \otimes_{\Z/n} G$ where the first arrow is $\Z/n(1) \otimes_{\Z/n} G \to \Z/n(1) \otimes_{\Z/n} G^{\et} \overset{h}\to H$. Then $\beta(h,a)$ is independent of the choice of $n$ and is defined globally on $X$. 

By the construction, we have functorial isomorphisms
\medskip

\ref{3.6}.1  \; $\beta(h, a_1+a_2) \cong \beta(h, a_1) \;(+)\; \beta(h, a_2)$
\medskip

\ref{3.6}.2 \; $\beta(h_1+h_2, a) \cong \beta(h_1, a) \;(+)\; \beta(h_2,a)$
\medskip

\noindent
where $(+)$ denote the Baer sums. ($h, h_i\in \Hom(G^{\et}(1), H)$, $a, a_i\in P^{\gp}$.) It follows that the functor 
$$\Hom(G^{\et}(1), H) \times P^{\gp} \to {\frak E}xt_{\log}(G, H)\;\;;\;\; (h,a) \mapsto \beta (h,a)$$
is extended, uniquely up to canonical isomorphisms, to a functor $\beta$ on $\Hom(G^{\et}(1), H) \otimes P^{\gp}$. 
\end{para} 

\begin{para}\label{3.7} We now define a 
quasi-inverse of \ref{3.4}.1. First we define a functor 
$$\gamma: {\frak E}xt_{\log}(G, H) \to \Hom(G^{\et}(1), H) \otimes P^{\gp}.$$
Let $E$ be an object of ${\frak E}xt_{\log}(G, H)$.
  Let $\epsilon: X^{\log}_{\fl}\to X^{\cl}_{\fl}$ be the canonical morphism of sites. By applying $R\epsilon_*$ 
to the exact sequence $0\to H\to E \to G \to 0$ on $X^{\log}_{\fl}$, we obtain an exact sequence on $X^{\cl}_{\fl}$
$$0\to H \to E\to G \overset{\delta}\to R^1\epsilon_*H.$$
We have $$R^1\epsilon_*H\cong \varinjlim_n  \cH om(\Z/n(1), H) \otimes  {\bf G}_{m,\log}/{\bf G}_m$$
by \cite{K2} 4.1. Let $x$ be the closed point of $X$ endowed with the inverse image of the log structure of $X$. The pullback of $\delta$ to $x$ is rewritten as $G|_x\to \varinjlim_n \cH om(\Z/n(1),H|_x) \otimes P^{\gp}$, 
where $G|_x$ and $H|_x$ denote the pullbacks on $x^{\log}_{\fl}$. This gives $G^{\et}(1)|_x \to H^{\mult}|_x\otimes P^{\gp}$, where $H^{\mult}$ 
denotes the multiplicative part of $H$. Since $\cH om(G^{\et}(1), H^{\mult})$ is represented by a finite \'etale scheme over $X$, we have 
$\Hom(G^{\et}(1), H^{\mult})= \Hom(G^{\et}(1)|_x, H^{\mult}|_x)$. Hence we obtain an element $\gamma(E)$ of $\Hom(G^{\et}(1), H) \otimes P^{\gp}$.
\end{para}

\begin{lem}\label{3.8}
(1) The composite 
$$\Hom(G^{\et}(1), H) \otimes P^{\gp}\to {\frak E}xt_{\log}(G,H) \to \Hom(G^{\et}(1), 
H)\otimes P^{\gp}$$ is the identity map.

(2) An object $E$ of ${\frak E}xt_{\log}(G,H)$ comes from ${\frak E}xt_{\cl}(G,H)$ if and only if $\gamma(E)\in \Hom(G^{\et}(1), H) \otimes P^{\gp}$ is zero.

\end{lem}

\begin{pf} (1) follows from the constructions of these maps.

(2) follows from the fact that the object $E$ comes from ${\frak E}xt_{\cl}(G,H)$ if and only if the sequence $0\to H \to E \to G \to 0$ on $X^{\cl}_{\fl}$ is exact.

\end{pf}

\begin{para}\label{3.9} Now we have a functor 

\medskip

\ref{3.9}.1   ${\frak E}xt_{\log}(G,H) \to {\frak E}xt_{\cl}(G, H) \times (\Hom(G^{\et}(1),H) \otimes P^{\gp}\;;$

$G\mapsto (G \;(-) \;\beta\circ \gamma(G), \; \gamma(G))$.

\medskip
\noindent
Here $(-)$ is the Baer subtraction. $G\; (-) \; \beta\circ \gamma(G)$ belongs to ${\frak E}xt_{\cl}(G, H)$ by \ref{3.8} (2), for its $\gamma$ is zero by \ref{3.8} (1). 

It is easily seen that the functors \ref{3.4}.1 and \ref{3.9}.1 are the quasi-inverses of each other. Thus we have proved Thm. \ref{3.3}. 
\end{para}

\begin{para}\label{3.10} We deduce Thm. \ref{3.1} from Thm. \ref{3.3}. An object $G$ of $\cC$ is identified with an object of 
${\frak E}xt_{\log}(G^{\et}, G^{\circ})$ (note \ref{2.3} 
that $(\fin/X)_d$ is stable under extensions  in the category of sheaves of abelian groups on $X^{\log}_{\fl}$), and an object $(G, N)$ of 
$\cC'$ is identified with an object of 
${\frak E}xt_{\cl}(G^{\et}, G^{\circ})\times \Hom(G^{\et}(1), G^{\circ})\otimes P^{\gp}$. 
A homomorphism $G\to H$ in $\cC$ (resp. $\cC'$) corresponds to a triple $(h_1, h_2, h_3)$ where $h_1$ is a 
homomorphism $G^{\circ}\to H^{\circ}$, $h_2$ is a 
homomorphism $G^{\et}\to  H^{\et}$ and $h_3$ is an 
isomorphism in ${\frak E}xt_{\log}(G^{\et}, H^{\circ})$ (resp. 
${\frak E}xt_{\cl}(G^{\et}, H^{\circ}) \times \Hom(G^{\et}(1), H^{\circ})\otimes P^{\gp}$)
from the object defined by $G$ and $h_1$ (by pushout)  to the object defined by $H$ and $h_2$ (by pullback). Thus the equivalence in \ref{3.3} gives the equivalence in \ref{3.1}.

\end{para}

\section{Logarithmic $p$-divisible groups}

Fix a prime number $p$. 

\begin{para}\label{4.1} We define the category $(\text{$p$-div}/X)_c$ (resp. $(\text{$p$-div}/X)_d$, resp. $(\text{$p$-div}/X)_r$, resp. $(\text{$p$-div}/X)_e$, resp. $(\text{$p$-div}/X)_f$) to be the full subcategory of the category of sheaves of abelian groups on $X^{\log}_{\fl}$ consisting of objects satisfying the following conditions (i) - (iii). Let $G_n=\Ker(p^n: G\to G)$. 

\medskip

(i)  $G=\cup_n\;  G_n$.

\medskip

(ii) $p: G\to G$ is surjective.

\medskip

(iii) For each $n\geq 0$, $G_n$ belongs to $(\fin/X)_c$ (resp. $(\fin/X)_d$, resp. $(\fin/X)_r$, resp. $(\fin/X)_e$, resp. $(\fin/X)_f$).

\medskip

We have
$$ (\text{$p$-div}/X)_c\subset (\text{$p$-div}/X)_d \subset (\text{$p$-div}/X)_r\subset (\text{$p$-div}/X)_e \subset (\text{$p$-div}/X)_f.$$

We will prove in \S7 later that $(\text{$p$-div}/X)_e = (\text{$p$-div}/X)_f$ (note that $(\fin/X)_e$ and $(\fin/X)_f$ do not coincide in general as was shown in \ref{1.8}.3). 

The category $ (\text{$p$-div}/X)_c$ is identified with the category of $p$-divisible groups over the underlying scheme of $X$. 

If $p$ is invertible on $X$, the categories  $(\text{$p$-div}/X)_d$, $(\text{$p$-div}/X)_r$, $(\text{$p$-div}/X)_e$ and $(\text{$p$-div}/X)_f$  coincide by \ref{2.1}. 

\end{para}

\begin{lem}\label{4.2} Let $G$ be an object of $(\text{$p$-div}/X)_f$. If $G_1$ belongs to $(\fin/X)_d$ 
(resp.   $(\fin/X)_r$, resp. $(\fin/X)_e$), then $G$ belongs to $(\text{$p$-div}/X)_d$ (resp.  $(\text{$p$-div}/X)_r$, resp. $(\text{$p$-div}/X)_e$). 

\end{lem}

This follows from \ref{2.3}.

\begin{para}\label{4.3} Relation with abelian varieties.

Let $K$ be a complete discrete valuation field and let $A$ be an abelian variety over $K$ such that $A\otimes_K L$ has semi-stable reduction for some finite extension $L$ of $K$ which is (at worst) tamely ramified over $K$. Let $S$ be the Spec of the valuation ring of $K$ endowed with the canonical log structure. Then the $p$-divisible group $A\{p\}$ over $K$ associated to $A$ extends to an object of $(\text{$p$-div}/S)_e$. It extends to an object of $(\fin/S)_d$ if either one of the following conditions (i) (ii) is satisfied.

(i) $A$ has semi-stable reduction.

(ii) $p$ is invertible on $X$. 

\end{para}

This \ref{4.3} was stated in the original version of this paper without proof. A proof was given by Zhao (\cite{Zh2} Thm. 1.2). We give a proof here for the convenience of the reader, in which we use the theory of log abelian varieties (\cite{KKN2})  developed by 
 T. Kajiwara, C. Nakayama, and the author. A proof of this style is already given by Zhao (\cite{Zh1} \S3.1).

\begin{para}\label{4.4} Let $K$ and $S$ be as in  \ref{4.3} and let $A$ be an abelian variety over $K$ of semi-stable reduction. 
 Then $A$ extends uniquely to a log abelian variety over $S$. 
 
 This follows from \cite{KKN6} 4.4. (A  correction to the proof of \cite{KKN6} 4.4 is given in  \cite{KKN7} C.2.)  (A proof of this is also given in Zhao \cite{Zh0} in the case $A$ is of split multiplicative reduction.)

\end{para}

By \ref{4.4}, the case $A$ is of semi-stable reduction of \ref{4.3} follows from \ref{4.5} below.

\begin{prop}\label{4.5} 
 Let  $A$ be a log abelian variety over $X$. 

(1) For $n\geq 1$, $A[n]:= \Ker(n: A\to A)$ belongs to $(\fin/X)_d$.

(2)  Let $p$ be a prime number. Then $A[p^{\infty}]:= \cup_n \; A[p^n]$ belongs to $(\text{$p$-div}/X)_d$.

\end{prop}

\begin{pf} If we replace $(\cdot)_d$ by $(\cdot)_r$, this is proved in 
\cite{KKN4} 18.1. It is enough to prove that the Cartier dual of $A[n]$ belongs to $(\cdot)_r$. By \ref{2pt}, it is sufficient to prove that the pullback of the Cartier dual of $A[n]$ to each point $x$ of $X$ belongs to $(\fin/x)_r$. But this pullback is $(A|_x)^*[n]$, where $(A|_x)^*$ is the dual log abelian variety of the pullback $A|_x$ of $A$ (\cite{KKN2} 7.5) and hence we are reduced to \cite{KKN4} 18.1. 
\end{pf}

\begin{para}\label{4.6}  We complete the proof of \ref{4.3}. Take a tamely ramified finite Galois extension $L$ of $K$ such that $A\otimes_K L$ is of semi-stable reduction. 
 Let  $T$ be the Spec of the valuation ring of $L$ endowed with the canonical log structure. 
By \ref{4.4} and \ref{4.5}, $A[p^{\infty}]\otimes_K L$ extends canonically to an object of $(\text{$p$-div}/T)_d$.
 Since $T\times_S T\cong \Gal(L/K) \times T$ for the fiber product in the category of fs log schemes, we have the descent data and hence this object of $(\text{$p$-div}/T)_d$ descends to an object of $(\text{$p$-div}/S)_e$. 

If $p$ is invertible on $S$, the last object belongs to $(\text{$p$-div}/S)_d$ by \ref{2.1}.

\end{para}

\begin{para}\label{4.7} By \cite{Zh1}  Thm. 1.1, a log abelian variety over $X$ is a sheaf on $X^{\log}_{\fl}$. Consider the following question. 

Let $K$ and $S$ be as in \ref{4.3}, and let $A$ be an abelian variety over $K$. Is it always true that $A$ extends to a sheaf $\cA$ of abelian groups on $S^{\log}_{\fl}$ which is locally on $S^{\log}_{\fl}$ a log abelian variety?

That is, we are asking whether we can generalize \ref{4.3} without the assumption about tame extension. 

However this question has a negative answer. (Actually the author expected that the answer was affirmative as is mentioned at the end of \cite{Zh2}.) 

In fact, for any integer $n$ which is invertible on $S$, $\cA[n]$ belongs to $(\fin/X)_f$ and by the assumption on $n$, $\cA[n]$ is locally constant on $X^{\log}_{\fl}$. By \cite{K2} 10.2, $\cA[n]$ is locally constant on $X^{\log}_{\et}$. This shows that  
 the action of $\Gal(\bar K/K)$ on $A[n]$ is tamely ramified. Of course, this is not true in general.

\end{para}

\section{Dieudonn\'e modules of logarithmic $p$-divisible groups}

Let $p$ be a prime number and assume that $p$ is nilpotent on $X$. The aim of this \S5
is to define under a certain condition on $X$, the {\it Dieudonn\'e crystal associated to a log $p$-divisible group over $X$. 

\begin{para}\label{5.1} The logarithmic crystalline site $(X/\Z_p)^{\log}_{\crys}$ is defined as follows (cf. \cite{K1} \S5). An object of $(X/\Z_p)^{\log}_{\crys}$ is a $4$-ple $(T, U, i, \gamma)$ where $T$ is an fs log scheme on which $p$ is locally nilpotent, $U$ is an \'etale scheme over the underlying scheme of $X$ endowed with the inverse image $M_U$ of $M_X$, $i$ is a morphism of fs log schemes $U\to T$ such that $i^*M_T \overset{\cong}\to M_U$ and such that the underlying morphism of schemes of $i$ is a closed immersion, and $\gamma$ is a PD ($=$ divided power) structure on the ideal $\Ker(\cO_T\to \cO_U)$ which is compatible with the PD structure on $p\Z_p\subset \Z_p$.

Morphisms in $(X/\Z_p)^{\log}_{\crys}$ are defined in the evident way. The topology on $(X/\Z_p)^{\log}_{\crys}$
is defined by the classical \'etale topology on $T$. The structure sheaf $\cO_{X/\Z_p}$ on $(X/\Z_p)^{\log}_{\crys}$ 
is defined by $\cO_{X/\Z_p}(T, U, i, \gamma)=\Gamma(T, \cO_T)$.

\end{para}

\begin{para}\label{5.2} Consider the following condition \ref{5.2}.1 on an object 
$T=(T, U, i, \gamma)$ of $(X/\Z_p)^{\log}_{\crys}$.

\medskip

\ref{5.2}.1. For every fs log scheme $T'$ over $T$ which is log flat and of Kummer type over $T$, there exists a PD structure $\gamma'$ on $\Ker(\cO_{T'}\to \cO_{U'})$, where $U'=U\times_T T'$, which is compatible with $\gamma$.

\medskip

Note that if $\gamma'$ exists, it is unique. Note also that $\gamma'$ exists if the underlying scheme of $T'$ is flat over that of $T$ in the classical sense (1.8 in \cite{Me} \S3). An example which does not satisfy \ref{5.2}.1 is given in \ref{5.10}.

\end{para}

\begin{para}\label{5.3} For an fs log scheme $X$, let $\lff(X)_f$ be the category of $\cO_X$-modules on $X^{\log}_{\fl}$ which are locally free of finite rank. 

We denote by $\lff(X)_c$ (resp. $\lff(X)_e$) the full subcategory of $\lff(X)_f$ 
consisting of objects which are locally free for the classical Zariski (resp. log \'etale) topology on $X$. The category $\lff(X)_c$ is equivalent to the category of $\cO_X$-modules on the underlying scheme of 
$X$ which are locally free of finite rank (for the classical Zariski topology).
\end{para}

\begin{para}\label{5.4} Let $X$ be as in \ref{5.1} and let $G$ be an object of $(\text{$p$-div}/X)_f$. For an object $T$ of $(X/\Z_p)^{\log}_{\crys}$ satisfying \ref{5.2}.1, we will define an object ${\bf D}(G)_T$ of $\lff(T)_f$ (\ref{5.3}) called the Dieudonn\'e module of $G$. It will be shown that if $G$ belongs to $(\text{$p$-div}/X)_d$ (resp. $(\text{$p$-div}/X)_e$), then ${\bf D}(G)_T$ belongs to $\lff(T)_c$ (resp. $\lff(T)_e$). (Note we will prove in \S7 that $(\text{$p$-div}/X)_e= (\text{$p$-div}/X)_f$.) Unlike in the classical Dieudonn\'e theory, in general, we can not extend the definition of ${\bf D}(G)_T$ to all objects of $(X/\Z_p)^{\log}_{\crys}$ (at least at present).

Let $T=(T, U, i, \gamma)$ be an object of $(X/\Z_p)^{\log}_{\crys}$ satisfying the 
condition \ref{5.2}.1. Working locally, take an integer $n\geq 2$ such that $p^n\cO_T=0$. We define ${\bf D}(G)_T$ by using just $G_n= \Ker(p^n:G\to G)$. Take a covering $T'\to T$ in $T^{\log}_{\fl}$ such that the pullback of $G_n$ to $U'=U\times_T T'$ belongs to $(\fin/U')_c$. Then $G_n\times_X U'$ is a truncated $p$-divisible group of 
level $n$ in the classical sense  (\cite{Me} 1.2, \cite{BBM} 3.3.8). Hence we have the contra-variant Dieudonn\'e module of 
$G_n \times_X U'$, which is an object of $\lff(T')_c$ endowed with operators $f$ and $v$ (\cite{BBM} 3.3.10). Let the covariant Dieudonn\'e module 
${\bf D}(G_n)_{T'}$ be the $\cO_{T'}$-dual of this contra-variant Dieudonn\'e module. It is endowed with the linear map 
$f:\varphi^*{\bf D}(G_n)_{T'} \to {\bf D}(G_n)_{T'}$ defined to be the $\cO_{T'}$-dual of the above $v$  and 
with the linear map ${\bf D}(G_n)_{T'}\to \varphi^*{\bf D}(G_n)_{T'}$ 
defined to be the $\cO_{T'}$-dual of the above $f$, such that $fv=p$ and $vf=p$.  (Here $\varphi$ is the Frobenius.)
 Furthermore, if $p_i$ ($i=1,2$) denote the first and the second projections $T' \times_T T' \to T'$, the identification of the two pullbacks of $G\times_X U'$ on $U'\times_U U'$ induces an 
 isomorphism $\theta: p_1^*({\bf D}(G_n)_{T'})\cong p_2^*({\bf D}(G_n)_{T'})$ in $\lff(T'\times_T T')_c$ which satisfies the standard cocycle condition on $T'\times_T T'\times_T T'$. 
 Hence by a formal argument on sheaves, we get an object $D={\bf D}(G_n)_T$ of $\lff(T)_f$ with $f: \varphi^*D\to D$ and $v: D\to \varphi^*D$ which 
 induces on $T'$ the $\cO_{T'}$-module ${\bf D}(G_n)_{T'}$ with $f$ and $v$.
 We define ${\bf D}(G)_T={\bf D}(G_n)_T$ (with $f$ and $v$ such that $fv=p$ and $vf=p$).  It is independent of the choices of $n$ and $T'\to T$.

For example, ${\bf D}(\Q_p/\Z_p)_T=\cO_T$ with $f=p$ and $v=1$, and ${\bf D}(\Q_p/\Z_p(1))_T=\cO_T$ with $f=1$ and $v=p$. 
\end{para}

\begin{lem}\label{5.5} If $G$ belongs to $(\text{$p$-div}/X)_d$ (resp. $(\text{$p$-div}/X)_e$), then ${\bf D}(G)_T$ belongs to $\lff(X)_c$ (resp. $\lff(X)_e$).

\end{lem}

\begin{pf} The statement for $(\text{$p$-div}/X)_e$ follows from the statement for  $(\text{$p$-div}/X)_d$ by the log \'etale localization. 

We prove the statement for $(\text{$p$-div}/X)_d$. We may work classical \'etale locally on $T$. Take $n$ such that $p^n\cO_T=0$. By \ref{2.7} (3), classical \'etale locally on $X$, we have an exact sequence $0\to H\to G_n \to H'\to 0$ of sheaves of abelian groups on $X^{\log}_{\fl}$ such that $H$ and $H'$ belong to $(\fin/X)_c$ and are classical truncated $p$-divisible group of level $n$. We have an exact sequence 
$0\to {\bf D}(H)_T \to {\bf D}(G)_T \to {\bf D}(H')_T \to 0$ in $\lff(T)_f$ (the exactness is checked log flat locally). Since ${\bf D}(H)_T$ and ${\bf D}(H')_T$ belong to $\lff(T)_c$, we see by \cite{K2} 6.6 (cf. \ref{7.2} of this paper) that ${\bf D}(G)_T$ also belongs to $\lff(T)_c$. 
\end{pf}

\begin{lem}\label{5.6} If $T$ and $T'$ are objects of $(X/\Z_p)^{\log}_{\crys}$ satisfying the condition \ref{5.2}.1 and if we are given a morphism $T'\to T$ of $(X/\Z_p)^{\log}_{\crys}$, 
then ${\bf D}(G)_{T'}$ is canonically isomorphic to the pullback of ${\bf D}(G)_T$ on $T'$.

\end{lem}

This is clear.

\begin{para}\label{5.7} In Thm. \ref{6.3}, we will use the following narrower log crystalline sites $(X/\Z_p)^{\log,\nar}_{\crys}$ and $(X/\Z_p)^{\log,\nar}_{\crys, \nilp}$:

Let  $(X/\Z_p)^{\log,\nar}_{\crys}$ be the full subcategory of $(X/\Z_p)^{\log}_{\crys}$ consisting of objects $(T,U, i, \gamma)$ satisfying 
the following condition \ref{5.7}.1.

\medskip

\ref{5.7}.1.  If $s$ is a local section of $M_T$ which is sent to zero in $\cO_U$, then $s$ is sent to zero in $\cO_T$.

\medskip
Let $(X/\Z_p)^{\log}_{\crys,\nilp}$ be the full subcategory of $X^{\log}_{\crys}$ 
consisting of objects $(T, U, i, \gamma)$ such that $\gamma$ is a nilpotent PD structure (that is, for each local section $x$ of the PD ideal on an affine open set, there exists $N\geq 1$ such that $\gamma_n(x)=0$ for every $n\geq N$.) Define the full subcategory 
$(X/\Z_p)^{\log,\nar}_{\crys, \nilp}$ of $X^{\log}_{\crys}$ to be  $(X/\Z_p)^{\log,\nar}_{\crys}\cap (X/\Z_p)^{\log}_{\crys, \nilp}$.

We regard these full subcategories of $(X/\Z_p)^{\log}_{\crys}$ as sites with the topologies 
  induced from the topology of $(X/\Z_p)^{\log}_{\crys}$.

If $M_X\to \cO_X$ is injective, then we have
$$(X/\Z_p)^{\log, \nar}_{\crys}= (X/\Z_p)^{\log}_{\crys}, \quad (X/\Z_p)^{\log, \nar}_{\crys, \nilp}= (X/\Z_p)^{\log}_{\crys, \nilp}$$
because there is no local section of $M_U$ which is sent to $0$ in $\cO_U$. 
\end{para}

\begin{para}\label{5.8} Now let $k$ be a perfect field of characteristic $p>0$ and assume that $X$ is over $k$ and satisfies the following condition \ref{5.8}.1.

\medskip

\ref{5.8}.1. Classical \'etale locally on $X$, there exist a torsion free fs monoid $P$, a chart $P\to M_X$, and an ideal $L$ of $P$ such that $L$ is sent to $(0)$ in $\cO_X$ and such that the induced morphism of schemes $X\to \Spec(k[P]/(L))$ is smooth in the classical sense. 

\medskip

In this case, for an object $G$ of $(\text{$p$-div}/X)_d$ (resp. $(\text{$p$-div}/X)_e$, resp. $(\text{$p$-div}/X)_f$), 
we have an object ${\bf D}(G)_T$ of $\lff(T)_c$ (resp. $\lff(T)_e$, resp. $\lff(T)_f$) for every object 
$T$ of $(X/\Z_p)^{\log, \nar}_{\crys}$ by \ref{5.5} and by the following \ref{5.9}. Hence we 
have a crystal ${\bf D}(G)$ of objects of $\lff(\cdot)_c$ (resp. $\lff(\cdot)_e$, resp. $\lff(\cdot)_f$) on $(X/\Z_p)^{\log,\nar}_{\crys}$ endowed with $f$ and $v$ such that $fv=p$ and $vf=p$.

Note that if $X$ is log smooth over $k$, then $X$ satisfies \ref{5.8}.1 with $L$ the empty ideal of $P$, and we have
$$(X/\Z_p)^{\log, \nar}_{\crys}= (X/\Z_p)^{\log}_{\crys}, \quad (X/\Z_p)^{\log, \nar}_{\crys, \nilp}= (X/\Z_p)^{\log}_{\crys, \nilp}.$$

\end{para}

\begin{prop}\label{5.9} Let $k$, $P$ and $X$ be as in \ref{5.8} and let 
 $T$ be an object of $X^{\log,\nar}_{\crys}$. Then 
 $T$ satisfies the condition \ref{5.2}.1.

\end{prop}

\begin{pf} Working classical \'etale locally, we may assume that 
$$X=U=\Spec(A),\quad  T=\Spec(R), \quad T'=\Spec(R')$$ 
in the commutative diagram 
$$\begin{matrix}  W(k)[P]/(L) &\to& A &\to & R & \to & A/pA\\
 \downarrow &&\downarrow &&\downarrow &&\\
 W(k)[Q]/(LQ) &\to & A' & \to & R'\end{matrix}$$
having the following properties.

$P$ is a torsion free fs monoid, $L$ is an ideal of $P$, $A$ is smooth over $W(k)[P]/(L)$, $Q$ is an fs monoid, the map $W(k)[P]/(L)\to W(k)[Q]/(LQ)$ is induced by an injective homomorphism of fs monoids $P\to Q$ of Kummer type, 
 two squares in the diagram are pushouts in the category of commutative rings, 
   the composition $A\to R \to A/pA$ is the canonical projection,
  the log structures of $X$ and $T$ are induced by the canonical log structure of $\Spec(W(k)[P])$, and the log structure of $T'$ is induced by the log structure of $\Spec(W(k)[Q])$.

 Take a surjective homomorphism $A[t_{\la}\;; \; \la\in\La]\to R$ 
 over $A$, where $t_{\la}$ are indeterminates. By replacing $t_{\la}$ by 
 $t_{\la}-a_{\la}$ for $a_{\la}\in A$ such that the images of $t_{\la}$ and $a_{\la}$ in $A/pA$ coincide, we may assume that the images of $t_{\la}$ in $A/pA$ are $0$. 
 Let $A\langle t_{\la} \;;\;\la\in \La\rangle$ be the PD polynomial ring over $A$. The PD structure $\gamma$ on $\text{Ker}(R\to A/pA)$ induces a ring homomorphism $A\langle t_{\la} \;;\;\la\in \La\rangle\to R$ over $A$ via which
 the canonical PD structure $\tilde \gamma$ on $\text{Ker}(A\langle t_{\la} \;;\;\la\in \La\rangle\to A/pA)$ is compatible with $\gamma$. 
  Our task is to prove that the PD structure on $\text{Ker}(R\to A/pA)$ extends to a PD structure on $\text{Ker}(R'\to A'/pA')$. Let $J=\text{Ker}(A\langle t_{\la}\;; \; \la\in\La\rangle\to R)$,  
 $J'=\text{Ker}(A'\langle t_{\la}\;; \; \la\in\La\rangle\to R')$. It is sufficient to prove that the canonical PD structure $\tilde \gamma'$ on $\text{Ker}(A'\langle t_{\la}\;; \; \la\in\La\rangle\to A'/pA')$ satisfies 
 $\tilde \gamma'_n(J')\subset J'$ for all $n\geq 1$. But this follows from the fact that $\tilde \gamma_n(J)\subset J$ for all $n\geq 1$ and the fact that the ideal $J'$ is generated by $J$. 
 \end{pf}

\begin{para}\label{5.10} Here we give an example of an object $T=(T,U, i, \gamma)$ of $(X/\Z_p)^{\log}_{\crys}$ which does not satisfy \ref{5.2}.1 and in which the underlying scheme of $X$ is $\Spec(k)$ for a perfect field $k$ of characteristic $p>0$. 

Let $p$ be a prime number and let $r>p$ be an integer. 
Let $Q$ be the  free commutative monoid generated by the two generators $x_0,x_1$, and let  $P$ be the
 submonoid of $Q$ generated by $y_j=x_0^jx_1^{r-j}$ for $0\leq j\leq r$. Then $P$ is an fs monoid because $P^{\gp}=\{x_0^ix_1^j\;|\; i,j\in \Z, i+j\equiv 0\bmod r\}$ and $P= P^{\gp}\cap Q$. Let $T=\Spec(k[P]/J^2)$ where $J$ is the ideal generated by $y_j$ ($0\leq j\leq r$), let $X=U= \Spec(k[P]/J)=\Spec(k)$, and endow these schemes with the inverse images of the canonical log structure of 
$\Spec(k[P])$. Let $i:U\to T$ be the inclusion map. Let $T'=\Spec(k[Q]) \times_{\Spec(k[P])} T=\Spec(k[Q]/J^2k[Q])$, $U'=T'\times_T U$ and endow these schemes with the inverse images of the canonical log structure of $\Spec(k[Q])$. Then $T'\to T$ is finite and log flat of Kummer type. It is log \'etale if $r$ is not divisible by $p$. However, $T'\to T$ is not flat.

Let $V$ be the $\F_p$-linear space with the base $y_j$ ($0\leq j\leq r$) and let $N$ be an $\F_p$-linear map $V\to V$ such that $N^2=0$. Then we have a PD structure $\gamma$ on the ideal $I=J/J^2=\Ker(\cO_T\to \cO_U)$ (here we identify a sheaf on a one-point set with the set of its global sections) defined as follows. Note $I= k\otimes_{\F_p} V$. We define $\gamma_0(x)=1$, $\gamma_1(x)=x$, $\gamma_n(x)=0$ unless $n=0,1,p$ for $x\in I$, $\gamma_p(\sum_{j=0}^r a_iy_j)= a_i^pN(y_j)$ for $a_i\in k$. 

We show that if $N$ is given as $N(y_0)= y_1$ and $N(y_j)=0$ for $1\leq j\leq r$, then $\gamma$ does not extend to a PD structure on $I'=\Ker(\cO_{T'}\to \cO_{U'})$. In  fact, if $\gamma$ extends to a PD structure $\gamma'$ on $I'$, since $x_0y_0=x_1y_1$, $\gamma'_p(x_0y_0)=x_0^p\gamma_p(y_0)=x_0^py_1$ must coincide with $\gamma'_p(x_1y_1)= x_1^p\gamma_p(y_1)=0$. But since $r>p$, $x_0^py_1$ is not zero in $I'$. 

\end{para}

\section{Categorical equivalence}

{\rm Here we state a result on the equivalence Thm. \ref{6.3} between the category of log $p$-divisible groups and the category of Dieudonn\'e modules with log poles under a certain assumption. We will give the proof in a forthcoming paper \cite{K3}.}

\begin{para}\label{6.1} For a prime number $p$ and for $X$ over $\F_p=\Z/p\Z$, we define
$$\cD\cC^{\log}_i(X)_c$$
with $i=1$ (resp. $i=2$) to be the category of crystals $D$ of objects of $\lff(\cdot)_c$ on $(X/\Z_p)^{\log,\nar}_{\crys}$ (resp. $(X/\Z_p)^{\log,\nar}_{\crys, \nilp}$)  
endowed with linear operators $f: \varphi^*D\to D$ 
and $v: D\to \varphi^*D$ satisfying $fv=p$, $vf=p$. Here $\varphi$ denotes the Frobenius. By replacing $\lff(\cdot)_c$ by $\lff(\cdot)_e$, we define the category $\cD\cC_i^{\log}(X)_e$.
\end{para}

\begin{para}\label{6.2} Let $p$ be a prime number and 
let $k$ be a perfect field of characteristic $p$. Assume that $X$ is over $W_n(k)$ for some $n\geq 1$ satisfying the following condition \ref{6.2}.1:

\medskip

\ref{6.2}.1.  Classical \'etale locally on $X$, there exist a chart $P\to M_X$ with $P$ a torsion free fs monoid and an ideal $L$ of $P$ such that the image of $L$ under $M_X\to \cO_X$ is zero and such that the morphism $X \otimes \F_p \to \Spec(k[P]/(L))$ of schemes is smooth.  

\medskip

We regard $X$ as an object of $(X\otimes\F_p)^{\log,\nar}_{\crys}$ with the PD structure $\gamma$ on $\text{Ker}(\cO_X\to \cO_{X\otimes \F_p}=\cO_X/p)$ defined by $\gamma_n(pa)= p^{[n]}a^n$ ($a\in \cO_X$, $p^{[n]}=p^n/n!$). (This $\gamma_n$ is well defined as is easily seen.)  If $p\neq 2$, we regard $X$ as an object of 
$(X\otimes\F_p)^{\log,\nar}_{\crys,\nilp}$.

Let $i=1$ or $2$. In the case $i=2$, assume $p\neq 2$. We denote by 
$$\cD\cC^{\log}_{h,i}(X)_c$$
the category of objects $D$ of $\cD\cC^{\log}_i(X\otimes \F_p)_c$ which are endowed with an $\cO_X$-submodule $H$ of $D_X$ satisfying the following conditions (i), (ii).

(i) $H$ is locally a direct summand of $D_X$.

\medskip
(ii) In $(\varphi^*D)_{X \otimes \F_p}$, $\varphi^*(H/pH)$ coincides with the image of
$$v: D_{X \otimes \F_p} \to (\varphi^*D)_{X \otimes \F_p}.$$

\medskip

We define $\cD\cC^{\log}_{h,i}(X)_e$ in the similar way.

\end{para}

\begin{thm}\label{6.3} Let $k$ be a perfect field of characteristic $p\neq 0, 2$, and assume that $X$ is over $W_n(k)$ for some $n\geq 1$ and satisfies \ref{6.2}.1.

(1) We have equivalences of categories
$$(\text{$p$-div}/X)_d \overset{\simeq}\to \cD\cC^{\log}_{h,1}(X)_c 
\overset{\simeq}\to \cD\cC^{\log}_{h,2}(X)_c,$$
$$(\text{$p$-div}/X)_e \overset{\simeq}\to \cD\cC^{\log}_{h,1}(X)_e 
\overset{\simeq}\to \cD\cC^{\log}_{h,2}(X)_e.$$

(2) If $p\cO_X=0$ (that is, if $X=X\otimes \F_p$) and if the underlying scheme of $X$ is reduced, the equivalences in (1) hold without the subscript $h$. 

\end{thm}

{\rm Note that for a log smooth scheme $X$ over $k$,
 $X$ satisfies \ref{6.2}.1 and satisfies also the assumptions in (2). Note also that in this case, we have 
$$(X\otimes \F_p)^{\log,\nar}_{\crys}= X^{\log}_{\crys}, \quad (X \otimes \F_p)^{\log,\nar}_{\crys,\nilp}= X^{\log}_{\crys,\nilp}$$ as in \ref{5.8}.}

\begin{para}\label{6.4}

Here is an outline of the proof of Thm. \ref{6.3} (1). The proof uses the unpublished preprint \cite{BK}. 

We have proved in \cite{BK} the part of this equivalence for $(\text{$p$-div}/X)_c$ and the full subcategories of $\cD \cC^{\log}_{h,i}(X)_c$ (resp. $\cD\cC^{\log}_i(X)_c$) ($i=1,2$) consisting of objects without log pole. 
We can reduce the above equivalence to this result of \cite{BK} as follows.

 \ref{6.3} is reduced to its "complete local version'', 
 that is, the categorical equivalence for 
 $\Spec({\hat \cO_{X,x}})$, for closed points $x$ of $X$. Assume $X$ is such Spec. The statement for $(\text{$p$-div}/X)_e$ is reduced 
 (by log \'etale localization) to the statement for $(\text{$p$-div}/X)_d$.  We consider
$(\text{$p$-div}/X)_d$. By \S3, $(\text{$p$-div}/X)_d$ 
is equivalent to pairs $(G,N)$ where $G$ is an object of $(\text{$p$-div}/X)_c$ and $N$ is a homomorphism $G^{\et}(1)\to G^{\circ}\otimes P^{\gp}$ for $P$ as in \ref{3.1}. 
We can prove a similar statement for $\cD\cC^{\log}_{h,i}(X)_c$: $\cD\cC^{\log}_{h,i}(X)_c$ is equivalent to the category of pairs 
$(D,N)$ where $D$ is an object of $\cD\cC^{\log}_{h,i}(X)_c$ 
with no log pole and $N$ is a homomorphism $D^{\et}(1) \to 
D^{\circ}\otimes P^{\gp}$ ($D^{\et}(1)$ and $D^{\circ}$ are defined in the natural way). By this,  we are reduced to \cite{BK}.
\end{para}

\section{Complements}

\begin{thm}\label{7.1} Let  $p$ be a prime number. We have
$$(\text{$p$-div}/X)_e= (\text{$p$-div}/X)_f.$$
\end{thm}

The proof is given in \ref{7.5} after preparations \ref{7.2}--\ref{7.4}.

\begin{lem}\label{7.2} Let $0\to \cF' \to \cF \to \cF'' \to 0$ be an exact sequence in $\lff(X)_f$. Then $\cF$ belongs to $\lff(X)_c$ (resp. $\lff(X)_e$) if and only if $\cF'$ and $\cF''$ belong to $\lff(X)_c$ (resp. $\lff(X)_e$).

\end{lem}

\begin{pf} This is 
reduced to the case $X$ is the Spec of 
a Noetherian strict local ring. In this case, by \cite{K2} 6.2, we have the following: 
Every object of $\lff(X)_f$ is a direct sum $\oplus_{i=1}^n \cL(a_i)$ for some 
 $a_i\in \Gamma(X,  M_X^{\gp}/\cO_X^\times\otimes \Q)$, where 
$\cL(a_i)$ is the line bundle on $X^{\log}_{\fl}$ corresponding to the ${\bf G}_m$-torsor $\{b\in {\bf G}_{m,\log}\;|\; b \bmod {\bf G}_m=a\}$ on $X^{\log}_{\fl}$
(note that on $X^{\log}_{\fl}$, we have ${\bf G}_{m,\log}/{\bf G}_m= {\bf G}_{m,\log}/{\bf G}_m \otimes \Q$). 
Furthermore, $\oplus_{i=1}^n \cL(a_i)\cong \oplus_{i=1}^n \cL(b_i)$ if and only if $n=m$ and 
$({\bar a}_1, \dots, {\bar a}_n)\bmod {\frak S_n}= ({\bar b}_1, \dots, {\bar b}_n) \bmod {\frak S}_n$, 
where ${\bar c}$ for $c\in \Gamma(X, M_X^{\gp}/\cO_X^\times\otimes \Q)$ denotes the image of $c$ in  
$\Gamma(X, M^{\gp}_X/\cO_X^\times \otimes \Q/\Z$). This proves \ref{7.2}. 
\end{pf}

\begin{prop}\label{7.3} Let $p$ be a prime number, let 
$G$ be an object of $(\text{$p$-div}/X)_f$, and consider the object $\Lie(G)$ of $\lff(X)_f$. Then $G$ belongs to $(\text{$p$-div}/X)_r$ (resp. $(\text{$p$-div}/X)_e$) if and only if $\Lie(G)$  belongs to $\lff(X)_c$ (resp. $\lff(X)_e$).
\end{prop}

\begin{pf} It is sufficient to consider the statement for $(\text{$p$-div}/X)_r$. 
By \cite{K2} 6.2 and by \ref{2pt}, we may (and do) assume that the underlying scheme of $X$ is the Spec of a field.

If $G$ belongs to $(\text{$p$-div}/X)_r$, then $G^{\circ}=\cup_n G_n^{\circ}$ belongs to $(\text{$p$-div}/X)_c$ (\ref{2.7} (2)), and hence $\Lie(G)= \Lie(G^{\circ})$ belongs to $\lff(X)_c$. 

Conversely,  assume that $\Lie(G)$ belongs to $\lff(X)_c$. Let $\cO_{G^{\circ}}= \varprojlim_n \cO_{G_n^{\circ}}$. Then locally on $X^{\log}_{\fl}$, we have 
$\cO_{G^{\circ}}\cong \cO_X[[T_1, \dots, T_d]]$. Let $I$ be the kernel of the homomorphism 
$\cO_{G^{\circ}}\to \cO_X$ 
defined by the origin of $G^{\circ}$, which is given  locally on $X^{\log}_{\fl}$ by $T_i \mapsto 0$ for all $i$. 
 Then  $I/I^2\cong \cH om_{\cO_X}(\Lie(G),\cO_X)$. By $\Lie(G)\in \lff(X)_c$, we have $I/I^2\in \lff(X)_c$. 
 We have $\text{Sym}^r_{\cO_X}(I/I^2) \overset{\cong}\to I^r/I^{r+1}$ 
 for all $r\geq 0$. Hence $I^r/I^{r+1}$ belongs to $\lff(X)_c$ 
 for every $r\geq 0$. From this, by \ref{7.2}, 
 we have that $\cO_{G^{\circ}}/I^r$ belongs to $\lff(X)_c$ for all $r\geq 0$. For each $n\geq 1$, $\cO_{G^{\circ}_n}$ is a quotient $\cO_X$-module 
 of $\cO_{G^{\circ}}/I^r$ for some $r\geq 1$, 
 and hence by \ref{7.2}, 
 we have that $\cO_{G^{\circ}_n}$ belong to $\lff(X)_c$ for all $n\geq 1$. Hence by \ref{2OG}, $G^{\circ}_n$ belong to $(\fin/X)_c$ for all $n\geq 1$, 
 and hence by \ref{2.7} (2), $G_n$ belong to $(\fin/X)_r$ for all $n\geq 1$. Hence $G\in (\text{$p$-div}/X)_r$. 
\end{pf}

\begin{lem}\label{7.4} Let $p$ be a prime number, let $r\geq 0$ and let $S$ be a finite subset of $(\Q/\Z)^r$ such that the set $\{ps\;|\;s\in S\}$ coincides with $S$. Then the order of every element of $S$ is coprime to $p$. 

\end{lem}

\begin{pf}
Assume that the order of some element of $S$ is divisible by $p$ and let $s$ be an element of $S$ which has the largest $p$-power part of the order. Then $s=ps'$ for some $s'\in S$. But the $p$-power part of the order of $s'$ is strictly bigger than that of $s$. A contradiction.
\end{pf}

\begin{para}\label{7.5} We prove \ref{7.1}. 
We may assume that the underlying scheme of $X$ 
is the Spec of a strict local ring. By \ref{2pt},
 we may assume that the underlying scheme of $X$ is $\Spec(k)$ for a separably closed field $k$. If $p$ is not the characteristic 
 of $k$, we are done by \ref{2.1}. Assume $p$ is the characteristic of $k$. We show that we may assume $k$ is algebraically closed.
  In fact, for the algebraic 
  closure $\bar k$ of $k$, if $\bar X$ denotes $\Spec(\bar k)$ endowed with the inverse image of the log structure 
  of $X$, 
  then by \cite{K2} 6.2, 
  $\Lie(G)$ belongs to $\lff(X)_e$ if and only if $\Lie(G) \otimes_k \bar k$ belongs to $\lff(\bar X)_e$, 
  and hence by \ref{7.3}, $G$ belongs to $(\text{$p$-div}/X)_e$ if and only if the pullback of $G$ to $\bar X$ belongs to $(\text{$p$-div}/\bar X)_e$. 

Now assume that   
 $X=\Spec(k)$ as a scheme for an algebraically closed field $k$ of characteristic $p>0$. 
Then  
for  an fs monoid $P$ such that $P^\times =\{1\}$, $X=Y\otimes \F_p$ with  $Y=\Spec(W(k))$ which is endowed with the log structure associated to 
$P\to W(k)\;;\; 1 \mapsto 1, P\smallsetminus \{1\}\to 0$. Let $\varphi:Y\to Y$ be the Frobenius which is compatible with the $p$-th power map on $P$.  
Then, for an object $G$ of $(\text{$p$-div}/X)_f$, ${\bf D}(G)$ gives an object $D$ of 
$\lff(Y)_f$.

Let $K$ be the field of fractions of $W(k)$  and let $\bar K$ be the algebraic closure of $K$. Let $Z=\Spec(K)$ and let $\bar Z=\Spec(\bar K)$ which are endowed with the inverse images of the log structure of $Y$. Consider the commutative diagram
$$\begin{matrix} \bar Z& \to & Z & \overset{\varphi}\to& Z \\ && \downarrow && \downarrow\\
&& Y & \overset{\varphi}\to& Y.
\end{matrix}$$

By  \cite{K2} 6.2, $D$ and the pullback $D_{\bar Z}$ of $D$ to $\bar Z$ 
correspond to the same  element $(a_i)_{1\leq i\leq n}\bmod {\frak S}_n$ of ${\frak S}_n\bs (P^{\gp}\otimes (\Q/\Z))^n)$ ($n= \text{rank}(D)$, 
$a_i\in P^{\gp}\otimes (\Q/\Z)$). On the other hand, $\varphi^*D$ and its pullback to $\bar Z$ correspond to $(a_i^p)_{1\leq i \leq n}\bmod {\frak S}_n$. (Here we denote the group law of $P^{\gp} \otimes \Q/\Z$ multiplicatively.)
Since the pullback of  $D$ on 
$Y$ to $Z$ and that of $\varphi^*D$ on $Y$ to $Z$ 
are isomorphic, their pullbacks to $\bar Z$ are isomorphic. 
This shows that  $(a_i)_i\bmod {\frak S}_n$ and  
 $(a^p_i)_i\bmod {\frak S}_n$ coincide. Now we apply \ref{7.4} as follows. We identify $P^{\gp}\otimes \Q/\Z$ with  $(\Q/\Z)^m$, where $m$ is the rank of $P^{\gp}$, and $(P^{\gp}\otimes \Q/\Z)^n$ with $(\Q/\Z)^r$ with $r=mn$. Write $a_i\in (\Q/\Z)^m$ as $(a_{i,j})_{1\leq j\leq m}$ ($a_{i,j}\in \Q/\Z$) and let $S= \{a_{i,j}\;|\; 1\leq i\leq n, 1\leq j\leq m\}\subset \Q/\Z$. Then $S$ satisfies the assumption of \ref{7.4} and hence the orders of $a_{i,j}$ are coprime to $p$. Hence the orders of $a_i$ as 
 elements of the group $P^{\gp}\otimes \Q/\Z$ are coprime to $p$  for all $i$.  
Hence $D$ belongs to $\lff(Y)_e$. Since $\Lie(G)$ is a quotient object of $D/pD$ in $\lff(X)_f$, 
 $\Lie(G)$ is an object 
of $\lff(X)_e$ (\ref{7.2}).
Hence by \ref{7.3},  $G$ belongs to $(\text{$p$-div}/X)_e$.
 \end{para}

\noindent Kazuya Kato

\noindent 
Department of Mathematics
\\
University of Chicago
\\
5734 S.\ University Avenue
\\
Chicago, Illinois, 60637 \\
USA
\\
\noindent kkato@math.uchicago.edu
\end{document}